\documentclass[12pt]{article}

\usepackage{amssymb,amsmath}

\usepackage{graphicx}

\usepackage{url}
\usepackage{enumerate}

\numberwithin{equation}{section}

\newenvironment{Proof}{\removelastskip\par\medskip
\noindent{\em Proof.} \rm}{\penalty-20\null\hfill$\square$\par\medbreak}

 \allowdisplaybreaks

 \def\complex{{\mathord{\mathbb C}}}
 \def\real{{\mathord{\mathbb R}}}
 
 \def\inte{{\mathord{\mathbb N}}}
 
 \def\qu{{\mathord{\mathbb Z}}}

 \def\real{{\mathord{{\rm I\kern-3pt R}}}}        
 
 \def\inte{{\mathord{{\rm I\kern-3pt N}}}}
 \def\sZZ{{\rm Z\kern-.45em{}Z}}
 \def\sQQ{{\kern 0.27em \vrule height1.45ex width0.03em depth0em
           \kern-0.30em \rm Q}}
 \def\qu{{\mathchoice
         {\sQQ}
         {\sQQ}
   {\kern 0.225em \vrule height1.05ex width0.025em depth0em \kern-0.25em \rm Q}
   {\kern 0.180em \vrule height0.78ex width0.020em depth0em \kern-0.20em \rm Q}
         }}
 \def\sGG{{\kern 0.27em \vrule height1.45ex width0.03em depth0em
           \kern-0.30em \rm G}}
 \def\gg{{\mathchoice
         {\sGG}
         {\sGG}
   {\kern 0.225em \vrule height1.05ex width0.025em depth0em \kern-0.25em \rm G}
   {\kern 0.180em \vrule height0.78ex width0.020em depth0em \kern-0.20em \rm G}
         }}

 \newtheorem{prop}{Proposition}[section]
 \newtheorem{lemma}[prop]{Lemma}
 
 \newtheorem{corollary}[prop]{Corollary}

\def\E{\mathop{\hbox{\rm I\kern-0.20em E}}\nolimits}

 \newcounter{hyp}
\title{\huge 
On the integral representations of $|\Gamma (z)|^2$ and its Fourier transform 
} 

\author{
Nicolas Privault 
\\ 
\normalsize 
Division of Mathematical Sciences 
\\ 
\normalsize 
School of Physical and Mathematical Sciences 
\\ 
\normalsize 
Nanyang Technological University 
\\ 
\normalsize 
SPMS-MAS-05-43, 21 Nanyang Link 
\\ 
\normalsize 
Singapore 637371
}

\begin{document}

\maketitle

\vspace{-0.9cm}

\baselineskip0.6cm
 
\begin{abstract} 
 We derive integral representations in terms of the Macdonald functions for the square modulus $s\mapsto | \Gamma ( a + i s ) |^2$ of the Gamma function and its Fourier transform when $a<0$ and $a\not= -1,-2,\ldots $, generalizing known results in the case $a>0$. This representation is based on a renormalization argument using modified Bessel functions of the second kind, and it applies to the representation of the solutions of the Fokker-Planck equation. 
\end{abstract} 
 

\baselineskip0.7cm
 
\section{Introduction} 
 The Fokker-Planck type equation 
\begin{equation} 
\label{bhe} 
\left\{
\begin{array}{l} 
 \displaystyle 
 \frac{\partial U_p}{\partial t} (t,y) =
 \left( 
 y^2 
 \frac{\partial^2}{\partial y^2} 
 + 
 y 
 \frac{\partial }{\partial y} 
 - 
 y^2 
 - 
 p^2 
 \right) 
 U_p (t,y), 
 \qquad
 y,t>0, 
\\ 
\\ 
 U_p ( 0 , y ) = y^p, 
\end{array} 
\right. 
\end{equation} 
 originates from statistical physics \cite{Schenzle}; 
 it is also connected to the analysis 
 of exponential functionals of Brownian motion \cite{ComtetMonthusYor}, 
 and to applications in mathematical finance, cf. 
 e.g. 
 \cite{carrschroder}, \cite{pintoux}, and references therein. 
 The solution of \eqref{bhe} can be 
 written using the heat kernel of the operator 
 $y^2 - y^2 \partial^2/\partial y^2 - y \partial/\partial y$ (see e.g. \cite{yakubovich1}) 
 as 
\begin{equation} 
\label{djksaldsa} 
 U_p ( t , y ) 
 = 
 \frac{2}{\pi^2} 
 \int_0^\infty  u \sinh(\pi u) 
 K_{iu}
 ( y ) 
 e^{ - ( p^2 + u^2 ) t } 
 \int_0^\infty 
 x^{p-1} 
 K_{iu} ( x ) 
 dx du. 
\end{equation} 
 Using the classical Mellin integral representation 
\begin{equation} 
\label{ir} 
 \Gamma \left( z + i s \right) 
 \Gamma \left( z - i s \right) 
 = 
 4 
 \int_0^\infty 
 \left( \frac{x}{2} \right)^{2z} 
 K_{2is} ( x ) 
 \frac{dx}{x} 
, 
 \qquad \Re ( z ) >0, \quad 
 s \in \real, 
\end{equation} 
 in the complex parameter $z\in (0,\infty ) + i \real$, 
 cf. e.g. relation~(26), page 331 of \cite{erdelyitable}, 
 one can show from the Fubini theorem that \eqref{djksaldsa} 
 can be turned for all $p>0$ into the single integral representation 
$$ 
 U_p ( t , y ) 
 = 
 \frac{2^p}{2 \pi^2} 
 \int_0^\infty  u \sinh(\pi u) 
 e^{ - ( p^2 + u^2 ) t } \left| 
 \Gamma\left( \frac{p }{2} + i \frac{u}{2} \right) \right|^2 
 K_{iu} ( y ) 
 du 
, 
 \quad 
 y,t >0, 
$$ 
 which is more suitable for numerical implementations. 
 Here, 
\begin{equation}
\label{u01} 
 K_w ( y )
 =
 \int_0^\infty
 e^{- y \cosh x }
 \cosh ( w x )
 dx 
,
 \qquad
 y > 0, 
\end{equation} 
 is the modifed Bessel function of the second kind, 
 or the Macdonald function,  
 with parameter $w\in \complex$, cf. relation~(5) 
 page~181 of \cite{Watson}. 
\\ 
 
 In this paper we give an extension of 
 the Kontorovich-Lebedev transform \eqref{ir} 
 to all non-integer negative values of $\Re ( z )$. 
 In particular, 
 when $z\in \complex$ with $\Re (z) \in (-1,0)$ and $s\in \real$, 
 we show that \eqref{ir} extends as 
\begin{eqnarray} 
\nonumber 
 \Gamma \left( z + i s \right) 
 \Gamma \left( z - i s \right) 
 \! \! 
 & = & 
 \! \! 
 4 
 \int_0^\infty 
 \left( 
 \frac{x}{2} 
 \right)^{2z} 
 K_{2is} ( x ) 
 \left( 
 1 
 + 
 \left( \frac{2}{x} \right)^{2z} 
 \frac{ 4 z }{\Gamma ( 1 -2z ) } 
 K_{2z} ( x ) 
 \right) 
 \frac{dx}{x} 
\\ 
 & & 
\label{ext} 
\end{eqnarray} 
 cf. Proposition~\ref{c04} below for the general result, which reads 
\begin{equation} 
\label{asdfg.1} 
 \Gamma \left( z + i s \right) 
 \Gamma \left( z - i s \right) 
 = 
 2 
 \int_0^\infty 
 K_{2is} ( x) 
 \left( 
 1 
 + 
 \left( 
 \frac{2}{x} 
 \right)^{2z} 
 \sum_{k=0}^n 
 \frac{ 4 ( k+z ) }{k! \Gamma ( 1- k - 2z )! } 
 K_{2k+2z} (x) 
 \right) 
 \left( 
 \frac{2}{x} 
 \right)^{1-2z} 
 \! \! \! \! \! \! 
 dx 
\end{equation}  
 when $z\in \complex$ with 
 $\Re (z) \in (-n-1,-n)$ and $s\in \real$, 
 for any $n \in \inte  = \{0,1,2,\ldots \}$. 
\\ 

 As an application, the representation \eqref{asdfg.1} allows 
 us to solve \eqref{bhe} for $p \in (-2n-2,0]$, $n\in \inte$, 
 by the Fubini theorem, as 
\begin{eqnarray} 
\nonumber 
 U_p ( s , y ) 
 & = & 
 \frac{2}{\pi^2} 
 \int_0^\infty  u \sinh(\pi u) 
 K_{iu}
 ( y ) 
 e^{ - ( p^2 + u^2 ) s } 
 \int_0^\infty 
 x^p 
 K_{iu} ( x ) 
 \frac{dx}{x} du 
\\ 
 & = &
 \frac{2}{\pi^2} 
 \int_0^\infty  u \sinh(\pi u) 
 K_{iu}
 ( y ) 
 e^{ - ( p^2 + u^2 ) s } 
\\ 
\nonumber 
 & & 
 \times 
 \int_0^\infty 
 x^p 
 K_{iu} ( x ) 
 \left( 
 1 
 + \sum_{k=0}^n 
 \frac{  2^{p+1} ( p+2k ) }{k! \Gamma ( 1- p - k ) } 
 K_{-p-2k} ( x ) 
 \right) 
 \frac{dx}{x} du 
\\ 
\nonumber 
 & & 
 - e^{-p^2 s } 
 \sum_{k=0}^n 
 \frac{  2^{p+1} ( p+2k ) }{k! \Gamma (1 -p - k )! } 
 K_{-p-2k} 
 ( y ) 
\\ 
\nonumber 
 & = & 
 \frac{2^p}{2 \pi^2} 
 \int_0^\infty  u \sinh(\pi u) 
 e^{ - ( p^2 + u^2 ) s } \left| 
 \Gamma\left( \frac{p }{2} + i \frac{u}{2} \right) \right|^2 
 K_{iu} ( y ) 
 du 
\\ 
\nonumber 
 & & 
 - 
 \sum_{k=0}^n 
 e^{ 4 ( ( p+2k)^2 -p^2 ) s } 
 \frac{  2^{p+1} ( p+2k ) }{k! \Gamma ( 1-p - k )! } 
 K_{-p-2k} ( y ) 
, 
\end{eqnarray} 
 $z\in \real$, $s >0$. 
 The above expression has been originally obtained in \cite{Schenzle} 
 using spectral expansions, 
 however the derivation presented here is much simpler 
 since the argument of \cite{Schenzle} 
 involves severe analytical difficulties in the computation of 
 normalization constants via the use of Meijer functions, 
 cf. page 1641 therein. 
\\ 
 
 On the other hand, Ramanujan showed in 
 \cite{ramanujan} that for $a>0$ the Fourier transform of 
 $s\mapsto | \Gamma ( a + i s )|^2$ 
 satisfies the relation 
\begin{equation} 
\label{rama} 
 \int_{-\infty}^\infty 
 e^{-i \xi s} 
 \left| 
 \Gamma \left( a + i s \right) 
 \right|^2 
 ds 
 = 
 \sqrt{\pi} 
 \Gamma (a) \Gamma ( a + 1/2) 
 ( \cosh ( \xi /2 ) )^{-2a} 
. 
\end{equation} 
 This relation has been extended to all $a\in (-1,0)$ 
 as an integral expression in Theorem~1.2 of \cite{chakrabarti}. 
\\ 
 
 As another application of \eqref{ext}, 
 we show that it can be used to deduce an extension of 
 \eqref{rama} to all non-integer 
 negative values of $a$, using integral expressions. 
 Namely when $a\in (-1,0)$ we  deduce the integral representation 
\begin{eqnarray} 
\label{djkld1} 
\lefteqn{ 
 \int_{-\infty}^\infty 
 e^{-i \xi s} 
 \left| 
 \Gamma \left( a + i s \right) 
 \right|^2 
 ds 
} 
\\ 
\nonumber 
 & = & 
 \frac{2\pi}{2^{2a}} 
 \int_0^\infty 
 x^{2a-1} 
 \left( 
 1 
 + 
 \left( \frac{2}{x} \right)^{2a} 
 \frac{ 4 a }{\Gamma ( 1 -2a ) } 
 K_{2a} ( x ) 
 \right) 
 e^{- x \cosh ( \xi/2 )} 
 dx, \qquad 
 \xi \in \real, 
\end{eqnarray} 
 for the Fourier transform of $s\mapsto | \Gamma ( a + i s )|^2$, 
 which is another extension of \eqref{rama} to $a\in (-1,0)$, 
 cf. Proposition~\ref{djkld1111} below for the general result 
 which reads 
\begin{eqnarray} 
\label{djkld1.1.2} 
\lefteqn{ 
 \int_{-\infty}^\infty 
 e^{-i \xi s} 
 \left| 
 \Gamma \left( a + i s \right) 
 \right|^2 
 ds 
} 
\\ 
\nonumber 
 & = & 
 \frac{2\pi}{2^{2a}} 
 \int_0^\infty 
 x^{2a-1} 
 \left( 
 1 
 + 
 \left( \frac{2}{x} \right)^{2a} 
 \sum_{k=0}^n 
 \frac{ 4 ( a + k ) }{k! \Gamma ( 1-k -2a ) } 
 K_{2a+2k} ( x ) 
 \right) 
 e^{- x \cosh (\xi/2 ) } 
 dx, 
\end{eqnarray} 
 $\xi \in \real$, $a\in (-n-1,-n)$, $n\in \inte$. 
\\ 
 
 The integrability as $x \to 0$ in \eqref{ext} and \eqref{djkld1} 
 is justified by the estimate 
\begin{equation} 
\label{jestimate} 
 x^{2a}  
 + 
 \frac{ 2^{2a+2} a }{\Gamma ( 1 -2a ) } 
 K_{2a} ( x ) 
 = 
 o ( x^{\varepsilon } ) 
, 
 \qquad 
 x \to 0, 
\end{equation} 
 for any $\varepsilon \in (0,2+2a) \cap (0,-2a)$, 
 when $a\in (-1,0)$,
 cf. \eqref{psi} and \eqref{irp}-\eqref{irp1} 
 below for the general case 
 $a \in (-n-1,-n)$, $n \in \inte$. 
\\ 
 
 This paper is organized as follows. 
 In Section~\ref{s2} we start by proving some asymptotic expansion 
 and integrability results that are needed for the proof of 
 both \eqref{ext} and \eqref{djkld1}. 
 The proof of the integral representation \eqref{asdfg.1} and its extension 
 to $\Re (z) \in (-n-1,-n)$ for all $n\in \inte$ are given 
 in Section~\ref{s3}. 
 Finally in Section~\ref{m} we derive the extension of 
 the Fourier transform identity 
 \eqref{djkld1} to $ \Re (z) \in (-n-1,-n)$ for all 
 $n\in \inte$ as a consequence of Proposition~\ref{c04}. 
\section{Asymptotic expansion and integrability} 
\label{s2} 
 In this section we derive 
 the asymptotic results needed for the proofs of 
 \eqref{ext}-\eqref{jestimate} in 
 Propositions~\ref{c04} and \ref{djkld1111} below. 
 We will use the modified Bessel function of the first kind 
\begin{equation} 
\label{1a} 
 I_z (x) = \sum_{k=0}^\infty 
 \frac{1}{\Gamma ( k + 1 ) \Gamma ( k+z + 1 ) } 
 \left( \frac{x}{2} \right)^{z+2k}, 
 \qquad 
 x \in \real, 
 \quad 
 z \in \complex. 
\end{equation} 
\begin{lemma} 
\label{l11.0} 
 For all $n\in \inte$,  $x\in \real$ and $z \in \complex$ we have 
\begin{eqnarray} 
\label{asd} 
\lefteqn{ 
 \left( 
 \frac{2}{x} 
 \right)^{2z} 
 \sum_{k=0}^n 
 \frac{k+z}{k! \Gamma (1- k - 2z )} 
 I_{2k+2z} (x) 
 = 
 \frac{\sin ( 2 \pi z ) }{2\pi} 
} 
\\ 
\nonumber 
 & & 
 \ \ \ \ \  \ \ \ \ \ \ \ \ \ \ \ 
 + 
 \sum_{l= n+1}^\infty 
 \frac{1}{l!} 
 \left( 
 \frac{x}{2} 
 \right)^{2l} 
 \sum_{k=0}^n 
 {l \choose k} 
 \frac{ k+z }{ \Gamma ( 1-2 z - k )! \Gamma ( k + l +2 z +1 )}. 
\end{eqnarray} 
\end{lemma} 
\begin{Proof} 
 From \eqref{1a} we have 
\begin{eqnarray} 
\nonumber 
\lefteqn{ 
 \! \! \! \! \! \! \! \! \! \! \! \! \! \! \! \! \! \! \! \! \! \! \! \! \! \! \! \! \! 
 \left( 
 \frac{2}{x} 
 \right)^{2z} 
 \sum_{k=0}^n 
 \frac{k+z}{k! \Gamma (1- k - 2z )} 
 I_{2k+2z} (x) 
 = 
 \sum_{k=0}^n 
 \frac{k+z}{k! \Gamma (1- k - 2z )} 
 \sum_{l=0}^\infty 
 \frac{1}{l! \Gamma ( 2k+l+2z+1) } 
 \left( \frac{x}{2} \right)^{2k+2l} 
} 
\\ 
\nonumber 
 & = & 
 \sum_{k=0}^n 
 \frac{k+z}{k! \Gamma (1- k- 2z )} 
 \sum_{l=k}^\infty 
 \frac{1}{(l-k)! \Gamma ( k+l+2z+1) } 
 \left( \frac{x}{2} \right)^{2l} 
\\ 
\label{dkldd} 
 & = & 
 \sum_{l=0}^\infty 
 \left( \frac{x}{2} \right)^{2l} 
 \sum_{k=0}^{\min ( n , l)}  
 \frac{k+z}{k! \Gamma (1- k - 2z ) (l-k)! \Gamma ( k+l+2z+1) } 
. 
\end{eqnarray} 
 Next, using Euler's reflection formula 
\begin{equation} 
\label{eul} 
 \displaystyle 
 \frac{1}{ \Gamma (k+l+2z+1) } 
 = 
 - (-1)^{k+l} 
 \frac{\sin ( 2 \pi z) }{\pi}  
 \Gamma(-2z-k-l), \qquad 
 k,l \in \inte, 
\end{equation} 
 cf. e.g. relation~(6.1.17), page 256 of \cite{Hand}, 
 we get 
\begin{eqnarray} 
\nonumber 
\lefteqn{ 
 \sum_{l=0}^n 
 \left( \frac{x}{2} \right)^{2 l } 
 \sum_{k=0}^l 
 \frac{k+z}{k!(l-k)! \Gamma (1- k - 2z ) \Gamma ( k + l +2 z +1 )} 
} 
\\ 
\nonumber 
 & = & 
 \frac{\sin ( 2\pi z) }{\pi} 
 \sum_{l=0}^n 
 \left( \frac{x}{2} \right)^{2 l } 
 \frac{(-1)^l}{l!} 
 \sum_{k=0}^l 
 (-1)^k 
 {l \choose k} 
 ( k+z ) 
 \frac{\Gamma ( -2 z - k - l )}{\Gamma (1-2z-k)} 
\\ 
\nonumber 
 & = & 
 \frac{\sin ( 2\pi z) }{\pi} 
 \sum_{l=0}^n 
 \left( \frac{x}{2} \right)^{2 l } 
 \frac{(-1)^l}{l!} 
 \frac{\Gamma ( -2 z - 2 l )}{\Gamma (1-2z-k)} 
 \sum_{k=0}^l 
 {l \choose k} 
 ( k+z ) 
 (-2z-2l)_{l-k} 
 (2z)_k 
, 
\\ 
\label{eqw1} 
\end{eqnarray} 
 where 
\begin{equation} 
\label{pk} 
 (p)_k = p (p+1) \cdots (p+k-1) 
\end{equation} 
 is the shifted factorial, $p\in \complex$, $k\geq 1$, with $(p)_0=1$. 
 Next, for all $l\geq 1$ and $z \in \complex$ we check that 
\begin{eqnarray} 
\nonumber 
\lefteqn{ 
 \sum_{k=0}^l 
 {l \choose k} 
 ( k+z ) 
 (-2z-2l)_{l-k} 
 (2z)_k 
} 
\\ 
\nonumber 
 & = & 
 z 
 \sum_{k=0}^l 
 {l \choose k} 
 (-2z-2l)_{l-k} 
 (2z)_k 
 - 
 k 
 \sum_{k=0}^l 
 {l \choose k} 
 (-2z-2l)_{l-k} 
 (2z)_k 
\\ 
\nonumber 
 & = & 
 z 
 \sum_{k=0}^l 
 {l \choose k} 
 (-2z-2l)_{l-k} 
 (2z)_k 
 - 
 \sum_{k=1}^l 
 \frac{l!}{(l-k)!(k-1)!} 
 (-2z-2l)_{l-k} 
 (2z)_k 
\\ 
\nonumber 
 & = & 
 z 
 \sum_{k=0}^l 
 {l \choose k} 
 (-2z-2l)_{l-k} 
 (2z)_k 
 + 
 2 z l 
 \sum_{k=0}^{l-1} 
 {l-1 \choose k} 
 (-2z-2l)_{l-1-k} 
 (2z+1)_{k}  
\\ 
\nonumber 
 & = & 
 z 
 (-2l)_l 
 + 
 2zl 
 (-2l+1)_{l-1} 
\\ 
\label{1} 
 & = & 0 
, 
\end{eqnarray} 
 where we used 
 the Pfaff-Saalsch\"utz 
 binomial identity 
$$ 
 (p+q)_l = 
 \sum_{k=0}^l {l \choose k} 
 (p)_k (q)_{l-k}, \qquad 
 p,q \in \complex, \quad l \in \inte,  
$$ 
 cf. e.g. Theorem~2.2.6 and Remark~2.2.1 of \cite{askey}. 
 As a consequence of \eqref{eqw1} and \eqref{1} we get 
$$ 
 \sum_{l=0}^n 
 \left( \frac{x}{2} \right)^{2 l } 
 \sum_{k=0}^l 
 \frac{k+z}{k!(l-k)! \Gamma (1- k - 2z ) \Gamma ( k + l +2 z +1)} 
 = 
 \frac{\sin ( 2\pi z) }{ 2 \pi} 
, 
$$ 
 which allows us to rewrite \eqref{dkldd} as 
\begin{eqnarray*} 
\lefteqn{ 
 \left( 
 \frac{2}{x} 
 \right)^{2z} 
 \sum_{k=0}^n 
 \frac{k+z}{k! \Gamma (1 - k -2z)} 
 I_{2k+2z} (x) 
} 
\\ 
 & = & 
 \frac{\sin ( 2\pi z) }{\pi} 
 + 
 \sum_{l=n+1}^\infty 
 \left( \frac{x}{2} \right)^{2 l } 
 \sum_{k= 0}^n 
 \frac{k+z}{k!(l-k)! \Gamma (1- k - 2z )! \Gamma ( k + 2 z + l +1)} 
, 
\end{eqnarray*} 
 and shows \eqref{asd}. 
\end{Proof} 
 As a consequence of Lemma~\ref{l11.0} we have the following 
 estimates. 
\begin{corollary} 
 Let $n\in \inte$. 
\begin{enumerate}[(i)] 
\item For all $z\in \complex$ such that $\Re ( z ) > -n-1$ we have 
\begin{equation} 
\label{pf} 
 \sum_{k=0}^n 
 \frac{k+z }{k! \Gamma ( 1- k - 2z )} 
 I_{2k+2z} (x) 
 = 
 \frac{\sin ( 2 \pi z )}{2 \pi} 
 \left( 
 \frac{x}{2} 
 \right)^{2z} 
 + 
 o ( x^{\varepsilon} ) 
, 
 \qquad 
 x \to 0, 
\end{equation} 
 for all $\varepsilon \in (0,2n+2+2\Re ( z ) )$. 
\item 
 For all $z\in \complex$ such that $-n-1 < \Re ( z ) <-n$ we have 
\begin{equation} 
\label{psi} 
 1 + \left( \frac{2}{x} \right)^{2z} 
 \sum_{k=0}^n 
 \frac{4(z+k) }{k! \Gamma ( 1- k - 2z )} 
 K_{2k+2z} (x) 
 = 
 o ( x^{\varepsilon - 2 \Re ( z)} ) 
, 
 \qquad 
 x \to 0, 
\end{equation}  
 for all $\varepsilon \in (0,2n+2+2 \Re (z) ) \cap (0,-2n-2 \Re (z) )$. 
\end{enumerate} 
\end{corollary} 
\begin{Proof} 
$(i)$ Relation~\eqref{pf} follows from Lemma~\ref{l11.0}. 
$(ii)$ On the other hand, relation~\eqref{1a} 
 with $-2 \Re (z) -2k\geq -2 \Re (z) -2n > \varepsilon > 0$ 
 shows that 
\begin{equation} 
\label{asfg} 
 I_{-2z-2k} (x) = 
 \sum_{l=0}^\infty 
 \frac{1}{l! \Gamma ( l -2z-2k +1 ) } 
 \left( \frac{x}{2} \right)^{{-2z-2k}+2l} 
 = 
 o ( x^\varepsilon ), 
 \qquad 
 x\to 0
, 
\end{equation} 
 hence \eqref{pf} and the identity 
\begin{equation} 
\label{id} 
 K_{2k+2z} (x) 
 = 
 \frac{\pi}{2 \sin ( 2\pi z )} 
 ( I_{-2z-2k}(x) - I_{2k+2z} (x) ), 
 \qquad 
 x \in \real, 
\end{equation} 
 allow us to conclude to \eqref{psi}. 
\end{Proof} 
 The next integrability result is a consequence of 
 Lemma~\ref{l11.0} and will be useful for the proofs of 
 Propositions~\ref{p05}, \ref{c04} and \ref{djkld1111} below. 
\begin{lemma} 
\label{p04} 
 Let $n\in \inte$. 
\begin{enumerate}[(i)] 
\item 
 For all $z\in \complex$ such that $-n-1 < \Re (z) $ 
 we have 
\begin{equation} 
\label{irp} 
 \sup_{s\in \real}
 \int_0^\infty 
 | K_{is}( x ) | 
 \left| 
 \frac{\sin ( 2 \pi z ) }{\pi} 
 -2 
 \left( \frac{2}{x} \right)^{2z} 
 \sum_{k=0}^n 
 \frac{ k+z }{k! \Gamma ( 1- 2 z - k )} 
 I_{2k+2z} (  x ) 
 \right| 
 \frac{d x }{ x^{1-2z}} < \infty
. 
\end{equation} 
\item 
 For all $z\in \complex$ such that $-n-1 < \Re (z) < -n$ we have 
\begin{equation} 
\label{irp1} 
 \sup_{s\in \real}
 \int_0^\infty 
 x^{2z-1} 
 | K_{is}( x ) | 
 \left| 
 1 
 + 
 \left( \frac{2}{x} \right)^{2z} 
 \sum_{k=0}^n 
 \frac{ 4 ( k+z ) }{k! \Gamma ( 1-k -2z ) } 
 K_{2k+2z} ( x ) 
 \right| 
 d x 
 < \infty
. 
\end{equation} 
\end{enumerate} 
\end{lemma} 
\begin{Proof} 
 $(i)$ 
 By relation~\eqref{u01}, 
 for all $\alpha >0$ there exists a constant $c_\alpha >0$ such that 
 $\cosh x > c_\alpha x^\alpha$ for all $x>0$, which shows that 
$$ 
 | K_{is}( y ) | 
 \leq 
 \int_0^\infty
 e^{- y \cosh x}
 dx
 \leq 
 \int_0^\infty
 e^{- y c_\alpha x^\alpha}
 dx 
 = 
 \frac{\Gamma\left(1/\alpha\right)}{\alpha \left(c_\alpha \right)^{1/\alpha}} 
 y^{-1/\alpha} 
, 
 \qquad 
 y>0, \quad s\in \real. 
$$ 
 Hence, using \eqref{pf} we have, 
 for all $s \in \real$ and $\alpha > 1/\varepsilon$, 
\begin{eqnarray*} 
\lefteqn{
 \! \! \! \! \! \! \! \! \! \! \! \! \! \! \! \! \! \! \! \! \! \! \! \! \! 
 \! \! \! \! \! \! \! \! \! \! \! \! \! \! \! \! 
 \int_0^1 | K_{is}( x ) | 
 \left| 
 \frac{\sin ( 2 \pi z ) }{\pi} 
 - 
 \left( 
 \frac{2}{x} 
 \right)^{2z} 
 \sum_{k=0}^n 
 \frac{2z+2k}{k! \Gamma ( 1- k - 2z ) } 
 I_{2k+2z} (x) 
 \right| 
 \frac{d x }{ x^{1-2z}} 
} 
\\ 
 & \leq & 
 c 
 \frac{\Gamma\left(1/\alpha\right)}{\alpha \left(c_\alpha \right)^{1/\alpha}} 
 \int_0^1
 \frac{d x }{ x^{1-\varepsilon+1/\alpha}}  
 < \infty, 
 \qquad 
 s\in \real, 
\end{eqnarray*} 
 for some constant $c>0$. 
 Next, the bound 
\begin{equation} 
\label{ntb1} 
 | K_{is} ( x ) | \leq K_0 ( x ), \qquad x>0, \quad s \in \real, 
\end{equation} 
 that follows from relation~\eqref{W} below and the equivalences 
\begin{equation} 
\label{ntb2} 
 K_{is} ( x ) \simeq e^{-x} \sqrt{\frac{\pi}{2x}},
 \qquad 
 x \to \infty, \quad s \in \real, 
\end{equation} 
 and 
\begin{equation} 
\label{ipx} 
 I_p ( x ) \simeq \frac{e^x}{\sqrt{2 \pi x}}, 
 \qquad x \to \infty, \quad p \in \complex, 
\end{equation} 
 show that 
$$ 
 \sup_{s\in \real} 
 \int_1^\infty 
 | K_{is}( x ) | 
 \left| 
 \frac{\sin ( 2 \pi z ) }{\pi} 
 -2 
 \left( \frac{2}{x} \right)^{2z} 
 \sum_{k=0}^n 
 \frac{ k+z }{k! \Gamma ( 1- 2 z - k )} 
 I_{2k+2z} (  x ) 
 \right| 
 \frac{d x }{ x^{1-2z}} < \infty
, 
$$ 
 which yields \eqref{irp} for all $z\in \complex$ such that 
 $\Re ( z) >-n-1$. 
\\ 
 
\noindent 
 $(ii)$ Due to the equivalences \eqref{id} and 
 \eqref{ipx} there also exists $c > 0$ such that 
\begin{equation} 
\label{psi2} 
 \left| 
 1 
 + 
 \left( \frac{2}{x} \right)^{2z} 
 \sum_{k=0}^n 
 \frac{ 4 ( k+z ) }{k! \Gamma ( 1-k -2z ) } 
 K_{2k+2z} ( x ) 
 \right| \leq c x^{-2z-1/2} e^x, 
\end{equation} 
 for sufficiently large $x>0$ and all $z\in \complex$, 
 and this yields \eqref{irp1} by replacing the use of \eqref{pf} 
 with that of \eqref{psi} in the proof of part $(i)$ above. 
\end{Proof} 
\section{Analytic continuation and integral representation} 
\label{s3} 
 In the next proposition, using analytic continuation, 
 we prove an integral representation formula 
 that will be applied to the proof of 
 \eqref{asdfg} in Proposition~\ref{c04} below. 
\begin{prop} 
\label{p05} 
 For all $n\in \inte$, $s\in \real$ 
 and $z\in \complex$ with $\Re (z) > -n-1$, $\Re (z) \notin -\inte$, 
 we have 
\begin{eqnarray} 
\label{a1.0} 
\lefteqn{ 
 \! \! \! \! \! \! \! 
 \Gamma \left( z + i s \right) 
 \Gamma \left( z - i s \right) 
 = 
 \frac{ 2 \pi  }{\sin ( 2 \pi z ) } 
 \sum_{k=0}^n 
 \frac{k+z}{k! ( (z+k)^2+s^2 ) \Gamma ( 1- k - 2z ) } 
} 
\\ 
\nonumber 
 & & 
 \! \! \! \! \! \! \! \! \! \! \! \! \! 
 \! \! \! \! \! 
 + 
 2^{2-2z} 
 \int_0^\infty 
 K_{2is} ( x) 
 \left( 
 1 
 - 
 \frac{ \pi}{\sin ( 2 \pi z ) } 
 \left( 
 \frac{2}{x} 
 \right)^{2z} 
 \sum_{k=0}^n 
 \frac{ 2z+2k }{k! \Gamma ( 1- k - 2z ) } 
 I_{2k+2z} (x) 
 \right) 
 \frac{dx}{x^{1-2z}} 
. 
\end{eqnarray} 
\end{prop} 
\begin{Proof} 
 Let $s\in \real \setminus \{0\}$. 
 We will prove the equality 
\begin{eqnarray} 
\label{djkldd1} 
\lefteqn{ 
 \! \! \! \! \! \! \! 
 \sin ( 2 \pi z ) 
 \Gamma \left( z + i s \right) 
 \overline{ \Gamma \left( \bar{z} + i s \right)} 
 = 
 \sin ( 2 \pi z ) 
 \Gamma \left( z + i s \right) 
 \Gamma \left( z - i s \right)
} 
\\ 
\nonumber 
 & = & 
 \pi 
 \sum_{k=0}^n 
 \frac{2z+2k}{k! ( (z+k)^2+s^2 ) \Gamma ( 1- k - 2z ) } 
\\ 
\nonumber 
 & & 
 \! \! \! \! \! \! \! \! \! \! \! \! \! 
 \! \! \! \! \! 
 + 
 \frac{4}{2^{2z}} 
 \int_0^\infty 
 K_{2is} ( x) 
 \left( 
 \sin ( 2 \pi z ) 
 - 
 \pi 
 \left( 
 \frac{2}{x} 
 \right)^{2z} 
 \sum_{k=0}^n 
 \frac{ 2z+2k }{k! \Gamma ( 1- k - 2z ) } 
 I_{2k+2z} (x) 
 \right) 
 \frac{dx}{x^{1-2z}}, 
\\ 
\label{djklaa} 
 & & 
\end{eqnarray} 
 for all $z = a+ i b \in \complex$ with 
 $a>-n-1$, in the following three steps. 
\\ 
 
\noindent 
$(i)$ Analyticity. 
 In \eqref{djkldd1}, the function 
$$ 
 \sin ( 2 \pi z ) 
 \Gamma \left( z + i s \right) 
 \Gamma \left( z - i s \right)
$$ 
 is analytic in $\{ z \ : \ z+is\notin ( -\inte ) , \ z-is \notin ( -\inte ) \}$ 
 and for each $k=0,1,\ldots ,n$ 
 the function 
 $( ( z +k)^2+s^2 )^{-1}$ is analytic in $z \in \complex \setminus (-\inte )$. 
 On the other hand, 
 by Lemma~\ref{l11.0} we can write the integrand 
 in \eqref{djkldd1} as 
\begin{eqnarray*} 
\lefteqn{ 
 x \mapsto 
 \frac{K_{is} ( x ) }{x^{1-2z}} 
 \left( 
 \frac{\sin ( 2 \pi z ) }{\pi} 
 - 
 \left( 
 \frac{2}{x} 
 \right)^{2z} 
 \sum_{k=0}^n 
 \frac{2z+2k}{k! \Gamma ( 1- k - 2z ) } 
 I_{2k+2z} (x) 
 \right) 
} 
\\ 
 & = & 
 - 
 2 K_{is} ( x ) 
 \sum_{k=0}^n 
 \frac{ x^{2z-1} }{k! \Gamma ( 1-2 z - k ) } 
 \sum_{l= n+1}^\infty 
 \left( 
 \frac{x}{2} 
 \right)^{2l} 
 \frac{ k+z }{ (l-k)! \Gamma ( k +2 z + l +1)} 
\\ 
 & = & 
 - 
 2 K_{is} ( x ) 
 \sum_{k=0}^n 
 \frac{ x^{2z-1} }{k! \Gamma ( 1-2 z - k )  \Gamma ( k +2 z )} 
 \sum_{l= n+1}^\infty 
 \left( 
 \frac{x}{2} 
 \right)^{2l} 
 \frac{ k+z }{ (l-k)! 
 ( k +2 z + l )
 \cdots 
 ( k +2 z )
} 
, 
\end{eqnarray*} 
 where 
 for each $k=0,1,\ldots ,n$ the partial 
 derivatives of 
$$ 
 z = a + ib \mapsto 
 \sum_{l=n+1}^\infty 
 \left( 
 \frac{x}{2} 
 \right)^{2l} 
 \frac{ k+z }{ (l-k)! 
 ( k + l +2 z )
 \cdots 
 ( k +2 z )
} 
$$ 
 with respect to $a$ and $b$ are locally 
 uniformly bounded by integrable functions of $x\in \real_+$ 
 from the bounds \eqref{ntb1} and \eqref{ntb2}, 
 by the same arguments as in the proof of Lemma~\ref{p04}. 
\\ 
 
 Hence we can exchange partial differentiation with respect to $a$ and 
 $b$ with the integration in \eqref{djklaa}, 
 showing that the Cauchy-Riemann conditions are satisfied by 
 the integral \eqref{djklaa} since all functions in the integrand 
 are analytic in $z\in \complex$. 
 Consequently, all terms in \eqref{djkldd1} are analytic in 
 $\{ z \in \complex \ : \ \Re (z ) > - n - 1 , \ 
 z+is\notin ( -\inte ) , \ z-is \notin ( -\inte ) \}$. 
\\ 
 
\noindent 
$(ii)$ The equality \eqref{djkldd1} holds 
 for all $s \in \real \setminus \{0\}$ 
 and  $z=a + i b \in (0,\infty ) + i \real$. 
 This follows from the integral representation \eqref{ir} which 
 reads 
$$ 
 \Gamma \left( z + i s \right) 
 \Gamma \left( z - i s \right) 
 = 
 4 
 \int_0^\infty 
 \left( \frac{x}{2} \right)^{2z} 
 K_{2is} ( x ) 
 \frac{dx}{x} 
, 
$$ 
 provided $a>0$, and from the Mellin transform 
\begin{equation} 
\label{djkldd} 
 4 \int_0^\infty 
 K_{2is} ( x) 
 I_{2k+2z} (x) 
 \frac{dx}{x} 
 = 
 \frac{ 1 }{(z+k)^2+s^2} 
\end{equation} 
 which is valid whenever $a +k > 0$, 
 cf. e.g. relation~(44) page 334 of \cite{erdelyitable}. 
\\ 
 
\noindent 
$(iii)$ 
 By analytic continuation the relation \eqref{djkldd1} extends 
 to $\{ z \in \complex \ : \ \Re (z ) > - n - 1 , \ 
 z+is \notin ( -\inte ) , \ z-is \notin ( -\inte ) \}$ 
 and we conclude by dividing \eqref {djklaa} by 
 $\sin ( 2 \pi z )$ when 
 $z\in \complex$ with $\Re (z) > -n-1$ and $\Re (z) \notin -\inte$. 
\end{Proof} 
 Note that in the above proof we could also have used 
 the unique continuation principle for real analytic 
 functions of $a$, see e.g. Corollary~1.2.3 of 
 \cite{krantz}, however real analyticity requires to check 
 the growth rate of partial derivatives, which would have been 
 more delicate. 
\\ 
 
 Relation~\eqref{ir} can be recovered from the integral representation 
\begin{equation} 
\label{W} 
 K_{is} ( y ) = \frac{1}{2} 
 \left( \frac{y}{2} \right)^{is}  
 \int_0^\infty 
 x^{-is-1} 
 e^{-x-y^2/(4x)} dx, \qquad 
 s,y \in \real, 
\end{equation} 
 cf. \cite{Watson} page~183, 
 using the Fubini theorem, as follows: 
\begin{eqnarray*} 
\nonumber 
 \Gamma \left( z + i s \right) 
 \Gamma \left( z - i s \right) 
 & = & 
 \int_0^\infty
 x^{-2is-1}e^{-x}
 x^{z+is}
 \int_0^\infty
 y^{-1+z+2is} 
 e^{-y}
 dy
 dx
\\
\nonumber 
 & = &
 \int_0^\infty
 \left( \frac{t}{2} \right)^{2z-1+2is}
 \int_0^\infty
 x^{-2is-1} e^{ -x-t^2/(4x) } 
 dx 
 dt
\\ 
\nonumber 
 & = & 
 4 
 \int_0^\infty 
 \left( \frac{t}{2} \right)^{2z} 
 K_{2is} ( t ) 
 \frac{dt}{t} 
, 
 \qquad s \in \real, 
\end{eqnarray*} 
 where we applied the change of variable $y = t^2 / (4x)$. 
 However, this argument is valid 
 only for $\Re (z) > 0$ due to integrability 
 restrictions in the exchange of integrals. 
\\ 
 
 We are now able to extend the above argument 
 to all $\Re (z) \in (-n-1,-n)$, $n \in \inte$, 
 in order to prove the integral representation 
 \eqref{asdfg} which also implies \eqref{ext}. 
\begin{prop} 
\label{c04} 
 For all $n\in \inte$, $s\in \real$ and $z\in \complex$ such that 
 $-n-1 < \Re (z ) < -n$ we have 
\begin{eqnarray} 
\label{asdfg} 
\lefteqn{ 
 \Gamma \left( z + i s \right) 
 \Gamma \left( z - i s \right) 
} 
\\ 
\nonumber 
 & = & 
 2 
 \int_0^\infty 
 K_{2is} ( x) 
 \left( 
 1 
 + 
 \left( 
 \frac{2}{x} 
 \right)^{2z} 
 \sum_{k=0}^n 
 \frac{ 4 ( k+z ) }{k! \Gamma ( 1- k - 2z )! } 
 K_{2k+2z} (x) 
 \right) 
 \left( 
 \frac{2}{x} 
 \right)^{1-2z} 
 \! \! \! \! \! \! 
 dx 
. 
\end{eqnarray} 
\end{prop} 
\begin{Proof} 
 First we note that the integrability in \eqref{asdfg} 
 follows from the bound \eqref{irp} above. 
 Next, from relations~\eqref{id}, \eqref{a1.0}, \eqref{djkldd} 
 with $-2\Re (z) -2k\geq -2 \Re ( z ) -2n > 0$, 
 we have 
\begin{eqnarray*} 
\nonumber 
\lefteqn{ 
 2^{2z-2} 
 \Gamma \left( z + i s \right) 
 \Gamma \left( z - i s \right) 
} 
\\ 
 & = & 
 \frac{\pi}{\sin ( 2 \pi z ) } 
 \int_0^\infty 
 K_{2is} ( x) 
 \left( 
 \frac{\sin ( 2 \pi z ) }{\pi} 
 - 
 \left( 
 \frac{2}{x} 
 \right)^{2z} 
 \sum_{k=0}^n 
 \frac{ 2z+2k }{k! \Gamma ( 1- k - 2z ) } 
 I_{2k+2z} (x) 
 \right) 
 \frac{dx}{x^{1-2z}} 
\\ 
 & & 
 + 
 \frac{2^{2z}\pi}{2\sin ( 2 \pi z )}  
 \sum_{k=0}^n 
 \frac{ k+z }{k! ( k+z )^2+s^2 ) \Gamma ( 1- k - 2z ) } 
\\ 
 & = & 
 \frac{\pi}{\sin ( 2 \pi z ) } 
 \int_0^\infty 
 K_{2is} ( x) 
 \left( 
 \frac{\sin ( 2 \pi z ) }{\pi} 
 - 
 \left( 
 \frac{2}{x} 
 \right)^{2z} 
 \sum_{k=0}^n 
 \frac{ 2z+2k }{k! \Gamma ( 1- k - 2z ) } 
 I_{2k+2z} (x) 
 \right) 
 \frac{dx}{x^{1-2z}} 
\\ 
 & & 
 + 
 \frac{2^z \pi}{\sin ( 2 \pi z ) } 
 \sum_{k=0}^n 
 \frac{ 2z+2k }{k! \Gamma ( 1- k - 2z ) } 
 \int_0^\infty 
 K_{2is} ( x) 
 I_{-2z-2k} (x) 
 \frac{dx}{x} 
\\ 
 & = & 
 \int_0^\infty 
 K_{2is} ( x) 
 \left( 
 1 
 + 
 \left( 
 \frac{2}{x} 
 \right)^{2z} 
 \sum_{k=0}^n 
 \frac{ 4 ( k+z ) }{k! \Gamma ( 1- k - 2z ) } 
 K_{2k+2z} (x) 
 \right) 
 \frac{dx}{x^{1-2z}} 
. 
\end{eqnarray*} 
\end{Proof} 
\section{Fourier transform of $| \Gamma ( a + i s ) |^2$} 
\label{m} 
 We begin by proving an integral representation for the Fourier transform 
 of 
$$ 
 s\mapsto | \Gamma ( a + i s )|^2, 
 \qquad 
 a\in (-n-1,-n), \quad n\in \inte, 
$$ 
 as a consequence of the integral representation 
 \eqref{asdfg.1} of Proposition~\ref{c04}. 
\begin{prop} 
\label{djkld1111} 
 Let $n\in \inte$ and $a\in (-n-1,-n)$. 
 For all $\xi \in \real$ we have 
\begin{eqnarray} 
\label{djkld1.1} 
\lefteqn{ 
 \int_{-\infty}^\infty 
 e^{-i \xi s} 
 \left| 
 \Gamma \left( a + i s \right) 
 \right|^2 
 ds 
} 
\\ 
\nonumber 
 & = & 
 \frac{2\pi}{2^{2a}} 
 \int_0^\infty 
 x^{2a-1} 
 \left( 
 1 
 + 
 \left( \frac{2}{x} \right)^{2a} 
 \sum_{k=0}^n 
 \frac{ 4 ( a + k ) }{k! \Gamma ( 1-k -2a ) } 
 K_{2a+2k} ( x ) 
 \right) 
 e^{- x \cosh (\xi/2 ) } 
 dx. 
\end{eqnarray} 
\end{prop} 
\begin{Proof} 
 This result can be informally deduced from \eqref{asdfg} 
 in Proposition~\ref{c04} 
 and the Fourier-Gelfand formula 
$$ 
 \int_{-\infty}^\infty 
 \cos ( 2 s y )
 e^{ -i \xi s} 
 ds 
 = 
 \pi 
 \left( 
 \delta ( \xi/2 - y ) 
 + 
 \delta ( \xi/2 + y ) 
 \right) 
$$ 
 in distribution theory, where 
 $\delta$ is the Dirac distribution at $0$. 
 However, with a view towards completeness, 
 we provide a proof by approximation 
 following the approach used in the proof of Theorem~1.1 of 
 \cite{chakrabarti}. 
 With the abbreviation 
\begin{equation} 
\label{psin} 
 \Psi_n ( x ) 
 : = 
 1 
 + 
 \left( \frac{2}{x} \right)^{2a} 
 \sum_{k=0}^n 
 \frac{ 4 ( a + k ) }{k! \Gamma ( 1-k -2a ) } 
 K_{2a+2k} ( x ), \qquad x\in \real, 
\end{equation} 
 we rewrite \eqref{asdfg} as 
$$ 
 \left| 
 \Gamma \left( a + i s \right) 
 \right|^2 
 = 
 \frac{4}{2^{2a}} 
 \int_0^\infty 
 x^{2a-1} 
 K_{2is} ( x ) 
 \Psi_n ( x ) 
 dx 
, 
 \qquad 
 s\in \real, 
$$ 
 for $a\in (-n-1,-n)$. 
 Then for any $\epsilon > 0$ we have 
\begin{eqnarray} 
\nonumber 
\lefteqn{ 
 \! \! \! \! \! \! \! \! \! \! \! \! \! \! \! \! \! \! \! \! \! 
 \int_{-\infty}^\infty 
 e^{-i \xi s} 
 \left| 
 \Gamma \left( a + i s \right) 
 \right|^2 
 ds 
 = 
 \int_{-\infty}^\infty 
 \lim_{\epsilon \to 0} 
 \left( 
 e^{- 2 \epsilon s^2 -i \xi s} 
 \left| 
 \Gamma \left( a + i s \right) 
 \right|^2 
 \right) 
 ds 
} 
\\ 
\nonumber 
 & = & 
 \lim_{\epsilon \to 0} 
 \int_{-\infty}^\infty 
 e^{- 2 \epsilon s^2 -i \xi s} 
 \left| 
 \Gamma \left( a + i s \right) 
 \right|^2 
 ds 
\\ 
\nonumber 
 & = & 
 2^{2-2a} 
 \lim_{\epsilon \to 0} 
 \int_{-\infty}^\infty 
 e^{- 2 \epsilon s^2 -i \xi s} 
 \int_0^\infty 
 x^{2a-1} 
 K_{2is} ( x ) 
 \Psi_n ( x ) 
 dx 
 ds 
\\ 
\label{dld} 
 & = & 
 2^{2-2a} 
 \lim_{\epsilon \to 0} 
 \int_0^\infty 
 x^{2a-1} 
 \Psi_n ( x ) 
 \int_{-\infty}^\infty 
 e^{- 2 \epsilon s^2 -i \xi s} 
 K_{2is} ( x ) 
 ds 
 dx 
, 
\end{eqnarray} 
 where the exchange of limit follows from the fact that 
 $s \mapsto  \left| \Gamma \left( a + i s \right) \right|^2 
 $ is a rapidly decreasing function in the Schwartz class, 
 and the last equality comes from \eqref{irp1} 
 below which ensures the integrability required for the 
 exchange of integrals. 
 Next, from relation~\eqref{u01} written as 
$$ 
 K_{2is} ( x ) 
 = 
 \int_0^\infty 
 e^{- x \cosh y } 
 \cos ( 2 s y ) 
 dy 
,
 \qquad
 x > 0, \quad 
 s \in \real
, 
$$ 
 we find 
\begin{eqnarray*} 
\lefteqn{ 
 \! \! \! \! \! \! \! \! \! \! \! \! \! \! \! 
 \int_{-\infty}^\infty 
 e^{- 2 \epsilon s^2 -i \xi s} 
 K_{2is} ( x ) 
 ds 
 = 
 \int_{-\infty}^\infty 
 e^{- 2 \epsilon s^2 -i \xi s} 
 \int_0^\infty
 e^{- x \cosh y }
 \cos ( 2 s y )
 dy 
 ds 
} 
\\ 
 & = & 
 \int_{-\infty}^\infty 
 e^{- x \cosh y }
 \int_0^\infty
 \cos ( 2 s y )
 e^{- 2 \epsilon s^2 -i \xi s} 
 ds 
 dy 
\\ 
 & = & 
 \frac{1}{4} \sqrt{\frac{\pi}{2\epsilon}} 
 \int_{-\infty}^\infty 
 e^{- x \cosh y } 
 ( 
 e^{-\frac{1}{2\epsilon } (y-\xi/2)^2 } 
 + 
 e^{
 -\frac{1}{2\epsilon } (y+\xi/2)^2 
 } 
 ) 
 dy 
, \qquad 
 x>0, 
\end{eqnarray*} 
 hence by \eqref{dld} we obtain 
\begin{eqnarray*} 
\lefteqn{ 
 \int_{-\infty}^\infty 
 e^{-i \xi s} 
 \left| 
 \Gamma \left( a + i s \right) 
 \right|^2 
 ds 
} 
\\ 
 & = & 
 2^{-2a} 
 \lim_{\epsilon \to 0} 
 \sqrt{\frac{\pi}{2\epsilon}} 
 \int_0^\infty 
 x^{2a-1} 
 \Psi_n ( x ) 
 \int_{-\infty}^\infty 
 e^{- x \cosh y } 
 ( 
 e^{
 -\frac{1}{2\epsilon } (y-\xi/2)^2 
 } 
 + 
 e^{ -\frac{1}{2\epsilon } (y+\xi/2)^2 
 } 
 ) 
 dy 
 dx 
\\ 
 & = & 
 2^{-2a} 
 \lim_{\epsilon \to 0} 
\sqrt{\frac{\pi}{2\epsilon}} 
 \int_0^\infty 
 x^{2a-1} 
 \Psi_n ( x ) 
 \int_{-\infty }^\infty 
 e^{- x \cosh (y+\xi/2 )  
 -\frac{1}{2\epsilon } y^2 
 } 
 dy 
 dx 
\\ 
 & & 
 + 
 2^{-2a} 
 \lim_{\epsilon \to 0} 
 \sqrt{\frac{\pi}{2\epsilon}} 
 \int_0^\infty 
 x^{2a-1} 
 \Psi_n ( x ) 
 \int_{-\infty}^\infty 
 e^{- x \cosh (y-\xi/2) 
 -\frac{1}{2\epsilon } y^2 
 } 
 dy 
 dx 
\\ 
 & = & 
 2^{-2a} 
 \sqrt{\frac{\pi}{2}} 
 \lim_{\epsilon \to 0} 
 \int_0^\infty 
 x^{2a-1} 
 \Psi_n ( x ) 
 \int_{-\infty}^\infty 
 e^{- x \cosh (z \sqrt{\epsilon} +\xi/2 )  
 -\frac{1}{2} z^2 
 } 
 dz 
 dx 
\\ 
 & & 
 + 
 2^{-2a} 
 \sqrt{\frac{\pi}{2}} 
 \lim_{\epsilon \to 0} 
 \int_0^\infty 
 x^{2a-1} 
 \Psi_n ( x ) 
 \int_{-\infty}^{\infty} 
 e^{- x \cosh (z\sqrt{\epsilon}-\xi/2) 
 -\frac{1}{2} z^2 
 } 
 dz 
 dx 
\\ 
 & = & 
 2^{1-2a} 
 \pi 
 \int_0^\infty 
 x^{2a-1} 
 \Psi_n ( x ) 
 e^{- x \cosh (\xi/2)} 
 dx, 
\end{eqnarray*} 
 where the required integrability follows from the bounds 
 \eqref{psi} and \eqref{psi2} of Section~\ref{s2}. 
\end{Proof} 
 In case $a>0$, $\Psi_{-1} ( x )$ in \eqref{psin} is 
 identically equal to $1$ and the proof of 
 Proposition~\ref{djkld1111} also yields the Mellin transform 
$$ 
 \int_{-\infty}^\infty 
 e^{-i \xi s} 
 \left| 
 \Gamma \left( a + i s \right) 
 \right|^2 
 ds 
 = 
 \frac{2\pi}{2^{2a}} 
 \int_0^\infty 
 x^{2a-1} 
 e^{- x \cosh ( \xi /2 ) } 
 dx
 = 
 \frac{2\pi}{2^{2a}} 
 ( \cosh ( \xi /2 ) )^{-2a} 
 \Gamma ( 2a ) 
, 
$$ 
 which recovers \eqref{rama}, 
 cf. also Theorem~1.1 in \cite{chakrabarti}. 
\\
 
 On the other hand, 
 when $a=-1/2$ the Fourier transform of $s \mapsto 
 | \Gamma ( -1/2 + i s )|^2$ can be explicitly computed as 
\begin{eqnarray*} 
 \int_{-\infty}^\infty 
 \! \! 
 e^{-i \xi s} 
 \left| 
 \Gamma \left( -1/2 + i s \right) 
 \right|^2 
 ds 
 & = & 
 4 \pi 
 \int_{-\infty}^\infty 
 \frac{ e^{-i \xi s}  }{(1+4s^2)\cosh ( \pi s )} ds 
\\ 
 & = & 
 4 \pi \log ( 1 + e^{-\xi} ) \cosh (\xi/2 ) + 
 2 \pi \xi e^{-\xi/2} 
, 
\end{eqnarray*} 
 cf. e.g. relation~(22) page 32 of \cite{erdelyitable}, 
 whereas \eqref{djkld1.1} yields 
$$ 
 \int_{-\infty}^\infty 
 e^{-i \xi s} 
 \left| 
 \Gamma \left( -1/2 + i s \right) 
 \right|^2 
 ds 
 = 
 4 
 \pi 
 \int_0^\infty 
 \left( 
 \frac{1}{x} 
 - K_{-1} ( x ) 
 \right) 
 e^{- x \cosh (\xi/2)} 
 \frac{dx}{x}, 
$$ 
 where the integrability in $0$ in the above integral 
 follows from 
$$ 
 \frac{1}{x} 
 - 
 K_{-1} ( x ) 
 = 
 o ( x^{\varepsilon} ) 
, 
 \qquad 
 x \to 0, 
$$ 
 for any $\varepsilon \in (0,1)$, cf. \eqref{jestimate} below. 

\footnotesize 

\def\cprime{$'$} \def\polhk#1{\setbox0=\hbox{#1}{\ooalign{\hidewidth
  \lower1.5ex\hbox{`}\hidewidth\crcr\unhbox0}}}
  \def\polhk#1{\setbox0=\hbox{#1}{\ooalign{\hidewidth
  \lower1.5ex\hbox{`}\hidewidth\crcr\unhbox0}}} \def\cprime{$'$}

\end{document}